\documentclass{amsart}

\newcommand{\R}{\ensuremath{ \mathbf{R} }}
\newtheorem{theorem}{Theorem}
\newtheorem{lemma}{Lemma}
\newtheorem{corollary}{Corollary}
\newcommand{\bt}{\begin{theorem}}
\newcommand{\et}{\end{theorem}}
\newcommand{\bl}{\begin{lemma}}
\newcommand{\el}{\end{lemma}}
\newcommand{\bc}{\begin{corollary}}
\newcommand{\ec}{\end{corollary}}
\newcommand{\beq}{\begin{equation}}
\newcommand{\eeq}{\end{equation}}
\newcommand{\benum}{\begin{enumerate}}
\newcommand{\eenum}{\end{enumerate}}
\newcommand{\card}{\ensuremath{\text{card}}}

\title[Relatively prime subsets of $\{m+1,\ldots,n\}$]{Asymptotic estimates for phi functions for subsets of $\{m+1, m+2,\ldots,n\}$}

\author{Melvyn B. Nathanson}
\address{Department of Mathematics\\
Lehman College (CUNY)\\
Bronx, New York 10468,
and School of Mathematics, Institute for Advanced Study,
Princeton, NJ 08540}
\email{melvyn.nathanson@lehman.cuny.edu, melvyn@ias.edu}
\thanks{The work of M.B.N. was supported in part by grants from the NSA Mathematical Sciences Program and the PSC-CUNY Research Award Program.}

\author{Brooke Orosz}
\address{Department of Mathematics\\ CUNY Graduate Center\\
New York, New York 10036}
\email{borosz@gc.cuny.edu}

\keywords{Relatively prime sets, Euler's phi function, Nathanson's phi function, combinatorial number theory, elementary number theory}
\subjclass[2000]{Primary 11A25, 11B05, 11B13, 11B75.} 

\date{\today}

\begin{document}

\begin{abstract}
Let $f(m,n)$ denote the number of  relatively prime subsets of $\{m+1,m+2,\ldots,n\}$, and let $\Phi(m,n)$ denote the number of subsets $A$ of $\{m+1,m+2,\ldots,n\}$ such that $\gcd(A)$ is relatively prime to $n$.  Let $f_k(m,n)$ and $\Phi_k(m,n)$ be the analogous counting functions restricted to sets of cardinality $k$.  Simple explicit formulas and asymptotic estimates are obtained for these four functions.
\end{abstract}

\maketitle

A nonempty set $A$ of integers is called \emph{relatively prime} if $\gcd(A) = 1.$    Let $f(n)$ denote the number of nonempty relatively prime subsets of $\{1,2,\ldots,n\}$ and, for $k \geq 1,$  let $f_k(n)$ denote the number of relatively prime subsets of $\{1,2,\ldots,n\}$ of cardinality $k$.  

Euler's phi function $\varphi(n)$ counts the number of positive integers $a$ in the set $\{1,2,\ldots,n\}$  such  that $a$ is relatively prime to $n$.   The Phi function $\Phi(n)$ counts the number of nonempty subsets $A$ of the set $\{1,\ldots,n\}$ such  that $\gcd(A)$ is relatively prime to $n$ or, equivalently, such that $A \cup\{n\}$ is relatively prime.  For every positive integer $k$, the function $\Phi_k(n)$ counts the number of sets $A \subseteq \{1,\ldots,n\}$ such  that $\card(A) = k$ and $\gcd(A)$ is relatively prime to $n$.

Nathanson~\cite{nath07a} introduced these four functions for subsets of $\{1,2,\ldots,n\}$, and El Bachraoui~\cite{elba07} generalized  them to subsets of the set $\{ m+1, m+2,\ldots, n\}$ for arbitrary nonnegative integers $m<n.$\footnote{Actually, our function $f(m,n)$ is El Bachraoui's function $f(m+1,n),$ and similarly for the other three functions.    This small change yields formulas that are more symmetric and pleasing esthetically.}   We shall obtain simple explicit formulas and asymptotic estimates for the four functions.

For every real number $x,$ we denote by $[x]$ the greatest integer not exceeding $x$.  We often use the elementary inequality $[x]-[y] \leq [x-y]+1$ for all $x,y \in \R.$

\bt   \label{RPS2:theorem:f}
For nonnegative integers $m < n,$ let $f(m,n)$ denote the number of relatively prime subsets of $\{m+1,m+2,\ldots,n\}.$  Then 
\[
f(m,n) = \sum_{d=1}^n \mu(d) \left( 2^{[n/d]-[m/d]} -1\right)
\]
and
\[
 0 \leq 2^{n-m} - 2^{[n/2]-[m/2]} - f(m,n)  \leq 2n2^{[(n-m)/3]}.
\]
\et

\begin{proof}
El Bachraoui~\cite{elba07} proved that
\[
f(m,n) = \sum_{d=1}^n \mu(d) \left( 2^{[n/d]} -1\right) - \sum_{i=1}^m\sum_{d|i}\mu(d)2^{[n/d]-i/d}.
\]
Rearranging this identity, we obtain 
\begin{align*}
f(m,n) & = \sum_{d=1}^n \mu(d) \left( 2^{[n/d]} -1\right) - \sum_{d=1}^m \mu(d) 2^{[n/d]} \sum_{\substack{i=1\\i|d}}^m2^{-i/d} \\
& = \sum_{d=1}^n \mu(d) \left( 2^{[n/d]} -1\right) - \sum_{d=1}^m \mu(d) 2^{[n/d]} \sum_{j=1}^{[m/d]} 2^{-j} \\
& = \sum_{d=1}^n \mu(d) 2^{[n/d]}  \left(1-\sum_{j=1}^{[m/d]} 2^{-j} \right)   
-\sum_{d=1}^n \mu(d) \\
& = \sum_{d=1}^n  \mu(d ) \left( 2^{[n/d]-[m/d]} -1 \right).
\end{align*}

Let $d \in \{1,2,\ldots,n\}.$  Then $m+1 \leq a \leq n$ and $d$ divides $a$ if and only if $[m/d]+1 \leq a/d \leq [n/d].$  It follows that $A\subseteq \{m+1,\ldots,n\}$ and $\gcd(A)=d$ if and only if $A' = (1/d)\ast A \subseteq \{[m/d]+1,\ldots,[n/d]\}$ and $\gcd(A')=1.$  Therefore,
\begin{align*}
2^{n-m}-1 
& = \sum_{d=1}^n f([m/d],[n/d]) \\
& \leq f(m,n) + 2^{[n/2]-[m/2]} - 1 + \sum_{d=3}^n 2^{[n/d]-[m/d]}
\end{align*}
and we obtain the lower bound
\[
f(m,n) \geq 2^{n-m}- 2^{[n/2]-[m/2]} - 2n 2^{ [(n-m)/3] }.
\]
For the upper bound, we observe that the number of subsets of even integers contained in the set $\{m+1,\ldots,n\}$ is exactly $2^{[n/2]-[m/2]}$ and so
\[
f(m,n) \leq 2^{n-m}- 2^{[n/2]-[m/2]}.
\]
This completes the proof.
\end{proof}

\bt      \label{RPS2:theorem:fk}
For nonnegative integers $m < n$ and for $k \geq 1$, let $f_k(m,n)$ denote the number of relatively prime subsets of $\{m+1,m+2,\ldots,n\}$ of cardinality $k$.  Then 
\[
f_k(m,n) = \sum_{d=1}^n \mu(d) {[n/d] - [m/d] \choose k}
\]
and
\[
0 \leq {n-m \choose k} - {[n/2]-[m/2] \choose k} - f_k(m,n) \leq n {[(n-m)/3] + 2\choose k}.
\]
\et

\begin{proof}
El Bachraoui~\cite{elba07} proved that
\[
f_k(m,n) = \sum_{d=1}^n \mu(d) {[n/d]  \choose k} - \sum_{i=1}^m\sum_{d|i}\mu(d) {[n/d] - i/d \choose k-1}.
\]
We recall the combinatorial fact that for $k \geq 1$ and $0 \leq M \leq N,$ we have 
\[
{N \choose k} - \sum_{j=1}^M {N-j \choose k-1} = {N-M \choose k}.
\]
Then 
\begin{align*}
f_k(m,n) & =  \sum_{d=1}^n \mu(d){[n/d]  \choose k} - \sum_{d=1}^m \mu(d) \sum_{\substack{i=1 \\ d|i }}^m {[n/d] - i/d \choose k-1}  \\
& =  \sum_{d=1}^m \mu(d)  \left(
 {[n/d]  \choose k} - \sum_{j=1}^{[m/d]} {[n/d] - j \choose k-1}
\right)  + \sum_{d=m+1}^n \mu(d){[n/d]  \choose k} \\
& =  \sum_{d=1}^m \mu(d) {[n/d] - [m/d] \choose k} + \sum_{d=m+1}^n \mu(d){[n/d]  \choose k} \\
& =  \sum_{d=1}^n \mu(d) {[n/d] - [m/d] \choose k}.
\end{align*}
We obtain an upper bound for $f_k(m,n)$ by deleting $k$-element sets of even integers:
\[
f_k(m,n) \leq {n-m \choose k} - {[n/2]-[m/2] \choose k} 
\]
and we obtain a lower bound from the identity 
\begin{align*}
{n-m \choose k} 
& = \sum_{d=1}^n f_k([m/d],[n/d]) \\
& \leq f_k(m,n) + {[n/2]-[m/2] \choose k} + \sum_{d=3}^n {[n/d]-[m/d] \choose k} \\
& \leq f_k(m,n) + {[n/2]-[m/2] \choose k} + n{[(n-m)/3] \choose k}.
\end{align*}
\end{proof}

\bt        \label{RPS2:theorem:Phi}
For nonnegative integers $m < n,$ let $\Phi(m,n)$ denote the number of subsets of $[m+1,n]$ such that $\gcd(A)$ is relatively prime to $n$.  Then
\[
\Phi(m,n) = \sum_{d|n} \mu(d)2^{(n/d)-[m/d]}.
\]
If ${p^{\ast}}$ is the smallest prime divisor of $n$, then
\[
0 \leq  2^{n-m} - 2^{(n/p^{\ast})-[m/p^{\ast}]} - \Phi(m,n)
\leq 2n 2^{ [ (n-m)/({p^{\ast}} +1) ] }.
\]
\et

\begin{proof}
El Bachraoui~\cite{elba07} proved that
\[
\Phi(m,n)  = \sum_{d|n} \mu(d) 2^{n/d} - \sum_{i=1}^m \sum_{d|(i,n)} \mu(d)2^{(n-i)/d}
\]
Rearranging this identity, we obtain
\begin{align*}
\Phi(m,n) & = \sum_{d|n} \mu(d) 2^{n/d} - \sum_{d|n} \mu(d) \sum_{\substack{i=1\\d|i}}^m 2^{(n-i)/d} \\
& = \sum_{d|n} \mu(d) 2^{n/d} - \sum_{d|n}\mu(d) \sum_{j=1}^{[m/d]} 2^{(n-jd)/d} \\
& = \sum_{d|n} \mu(d) 2^{n/d}\left[1 -  \sum_{j=1}^{[m/d]} 2^{-j}\right] \\
& = \sum_{d|n} \mu(d) 2^{(n/d)-[m/d]}.
\end{align*}
Let ${p^{\ast}}$ be the smallest prime divisor of $n$.  Deleting all subsets of $\{m+1,\ldots,n\}$ whose elements are all multiplies of $p^{\ast}$, we obtain the upper bound
\[
\Phi(m,n) \leq  2^{n-m} - 2^{(n/p^{\ast})-[m/p^{\ast}]}.
\]
For the lower bound, we have
\begin{align*}
\Phi(m,n)&  - \left( 2^{n-m} - 2^{(n/p^{\ast})-[m/p^{\ast}]}  \right)
 = \sum_{ \substack{ d|n \\ d >  {p^{\ast}}} } \mu(d) 2^{(n/d)-[m/d]}   \\
& \leq 2 \sum_{ \substack{ d|n \\ d >  {p^{\ast}}} } 2^{ [ (n-m)/d ] }   
\leq 2n 2^{ [ (n-m)/({p^{\ast}} +1) ] }.
\end{align*}
This completes the proof.
\end{proof}

\bt       \label{RPS2:theorem:Phik}
For nonnegative integers $m < n,$
let $\Phi_k(m,n)$ denote the number of subsets of cardinality $k$ contained in the interval of integers  $\{m+1,m+2,\cdots n\}$ such that $\gcd(A)$ is relatively prime to $n$.   Then
\[
\Phi_k(m,n) =  \sum_{d|n} \mu(d) {n/d - [m/d] \choose k} 
\]
and
\[
0 \leq {n-m  \choose k} -  {n/p^{\ast} - [m/p^{\ast}] \choose k} - \Phi_k(m,n)   \leq  n{[(n-m)/ (p^{\ast}+1) ] +1\choose k}.
\]
\et

\begin{proof}
Let $p^{\ast}$ be the smallest prime divisor of $n.$  
El Bachraoui~\cite{elba07} proved that
\[
\Phi_k(m,n) = \sum_{d|n} \mu(d) {n/d \choose k} - \sum_{i=1}^m \sum_{d|\gcd(i,n)} \mu(d) {(n-i)/d \choose k-1 }.
\]
Rearranging this identity, we obtain
\begin{align*}
\Phi_k(m,n) 
& = \sum_{d|n} \mu(d) {n/d \choose k} - \sum_{d|n} \mu(d) \sum_{\substack{ i=1 \\ i|d}}^m  {(n-i)/d \choose k-1 } \\
& = \sum_{d|n} \mu(d) \left(  {n/d \choose k} - \sum_{j=1}^{[m/d]} {n/d - j \choose k-1 } \right) \\
& = \sum_{d|n} \mu(d) { n/d - [m/d] \choose k} \\
& \geq  {n - m \choose k} -  {n/p^{\ast} - [m/p^{\ast}] \choose k} 
  - \sum_{\substack{d|n \\ d > p^{\ast} }} {n/d - [m/d] \choose k} \\
& \geq  {n - m \choose k} -  {n/p^{\ast} - [m/p^{\ast}] \choose k} 
  - \sum_{\substack{d|n \\ d > p^{\ast} }} { [(n-m)/d] +1\choose k} \\
& \geq  {n - m \choose k} -  {n/p^{\ast} - [m/p^{\ast}] \choose k} 
  - n{[(n-m)/(p^{\ast}+1)] +1\choose k}.
\end{align*}
Deleting $k$-element subsets of $\{m+1,\ldots,n\}$ whose elements are multiples of $p^{\ast}$, we get the upper bound
\[
\Phi_k(m,n) \leq {n-m \choose k} - {[n/p^{\ast}]-[m/p^{\ast}] \choose k}.
\]
This completes the proof.  
\end{proof}

\def\cprime{$'$} \def\cprime{$'$} \def\cprime{$'$}
\providecommand{\bysame}{\leavevmode\hbox to3em{\hrulefill}\thinspace}
\providecommand{\MR}{\relax\ifhmode\unskip\space\fi MR }
\providecommand{\MRhref}[2]{%
  \href{http://www.ams.org/mathscinet-getitem?mr=#1}{#2}
}
\providecommand{\href}[2]{#2}

\end{document}